\documentclass[12pt,a4paper]{amsart}
\usepackage{amssymb}
\usepackage{mathrsfs}
\headsep=1cm \numberwithin{equation}{section}
\newtheorem{thm}{Theorem}[section]
\newtheorem{cor}[thm]{Corollary}
\newtheorem{lem}[thm]{Lemma}
\newtheorem{prop}[thm]{Proposition}
\theoremstyle{definition}

\theoremstyle{remark}

\theoremstyle{example}
\newtheorem{exam}[thm]{Example}
\numberwithin{equation}{section}

\newcommand\Supp{\operatorname{Supp}}
\newcommand\Ass{\operatorname{Ass}}
\newcommand\mAss{\operatorname{mAss}}

\newcommand\Spec{\operatorname{Spec}}
\newcommand\Rad{\operatorname{Rad}}
\newcommand\Hom{\operatorname{Hom}}

\newcommand\grade{\operatorname{grade}}
\newcommand\height{\operatorname{ht}}

\newcommand\m{\operatorname{\frak m}}
\newcommand\q{\operatorname{\frak q}}
\newcommand\p{\operatorname{\frak p}}

\begin{document}
\title[Faltings'  finiteness dimension of local cohomology modules ]
{Faltings'  finiteness dimension of local cohomology modules over local Cohen-Macaulay rings}
\author{Kamal Bahmanpour and Reza Naghipour$^*$ }

\address{Department of Mathematics,  Faculty of Sciences, University of Mohaghegh Ardabili, 56199-11367,  Ardabil, Iran,
 and School of Mathematics, Institute for Research in Fundamental Sciences (IPM), P.O. Box: 19395-5746, Tehran, Iran.}
 \email{bahmanpour.k@gmail.com}
 \address{Department of Mathematics, University of Tabriz, Tabriz, Iran;
and School of Mathematics, Institute for Research in Fundamental Sciences (IPM), P.O. Box: 19395-5746, Tehran, Iran.}
\email{naghipour@ipm.ir} \email {naghipour@tabrizu.ac.ir}%
\thanks{ 2010 {\it Mathematics Subject Classification}: 13D45, 14B15.\\
This research was in part supported by a grant from IPM (No.  92130022).\\
$^*$Corresponding author: e-mail: {\it naghipour@ipm.ir} (Reza Naghipour)}%
\keywords{Cohen Macaulay ring, equidimensional ring, finiteness dimension, local cohomology}
\begin{abstract}
Let $(R, \frak m)$ denote a local Cohen-Macaulay ring and $I$ a non-nilpotent ideal of $R$. The purpose of this article is to investigate  Faltings' finiteness
dimension $f_I(R)$ and equidimensionalness of certain homomorphic image of $R$. As a consequence
we deduce that $f_I(R)={\rm max}\{1, \height I\}$ and if  $\mAss_R(R/I)$ is cotained  in  $\Ass_R(R)$, then the ring   $R/ I+\cup_{n\geq 1}(0:_RI^n)$ is  equidimensional of dimension $\dim R-1$.
Moreover, we will obtain a lower bound for injective dimension of the local cohomology module $H^{\height I}_I(R)$, in the case
$(R, \frak m)$ is a complete  equidimensional local ring.
\end{abstract}
\maketitle
\section{Introduction}
Throughout this paper, let $R$ denote a commutative Noetherian
ring (with identity) and $I$ an ideal of $R$. For an $R$-module $L$,
the $i^{th}$ local cohomology module of $L$ with respect to $I$ is
defined as$$H^i_I(L) = \underset{n\geq1} {\varinjlim}\,\,
\text{Ext}^i_R(R/I^n, L).$$ We refer the reader to \cite{Gr1} or \cite{BS} for more details about local cohomology.

For any finitely generated $R$-module $M$, the notion $f_I(M)$, the {\it finiteness dimension} of $M$ relative to $I$, is defined
to be the least integer $i$ such that $H^i_I(M)$ is not finitely generated, if there exist such $i$'s and $\infty$ otherwise, i.e.
$$f_I(M):=\inf\{i\in \Bbb{N}_0|\,H^i_I(M)\,\,{\rm is}\,\,{\rm not}\,\,{\rm finitely}\,\,{\rm generated} \}.$$

Our objective in this paper is to investigate the finiteness dimension $f_I(R)$, when $R$ is a local Cohen-Macaulay ring. More precisely, as
a main result we shall show that:

\begin{thm} Let $(R,\mathfrak{m})$ be a   Cohen-Macaulay local ring and  $I$ a  non-nilpotent ideal of $R$. Then
 $f_I(R)={\rm max}\{1,\height I \}.$
\end{thm}

One of our tools for proving Theorem 1.1 is the following, which will play a key role in the proof of that theorem.

\begin{prop}
Let $(R,\mathfrak{m})$ be a  Cohen-Macaulay local ring and let $X$ and $Y$ be non-empty subsets of $\Ass_R(R)$ such that $\Ass_R(R)=X\cup Y$ and $X\cap Y=\emptyset$.  Then  $R/(I+J)$ is an equidimensional local ring of dimension $\dim R-1$, where
 $I=\cap_{\p\in X}\p$ and $J=\cap_{\p \in Y}\p.$
\end{prop}

Recall that a Noetherian ring $R$, of finite Krull dimension $d$, is called {\it equidimensional} if $\dim R/\p=d$ for every minimal
prime ideal $\p$ of $R$. As an another main result, we shall show that:

\begin{thm}
Let $(R,\mathfrak{m})$ be a  Cohen-Macaulay local ring  and let $I$ be a non-nilpotent ideal of $R$ such that
$\mAss_R(R/I)\subseteq \Ass_R(R)$. Then  $R/ (I+\cup_{n\geq 1}(0:_RI^n))$ is an equidimensional local ring of dimension $\dim R-1$.
\end{thm}

In \cite{HaS}, Hartshorne and Speiser, proved  that if $(R, \m, k)$ is a regular local ring, contains a field of characteristic $p>0$, and $H^i_I(R)$
is supported only at the maximal ideal, then $\Hom_R(k, H^i_I(R)$ is  a finitely generated $R$-module and, moreover, $H^i_I(R)$ is injective.
Also, Huneke and Sharp in \cite{HS} made a remarkable breakthrough. They generalized Hartshorne-Speiser's result by proving that if $R$ is
any regular ring containing a field of characteristic $p>0$, then ${\rm inj}\dim H_I^i(R)\leq \dim\Supp H_I^i(R),$ where ${\rm inj}\dim H_I^i(R)$  denotes
the injective dimension of $H_I^i(R)$ and $\dim\Supp H_I^i(R)$ stands for the dimension of the support of $H_I^i(R)$ in $\Spec R$. Finally, in
\cite{Ly1}, Lyubeznik generalized the above-mentioned result of Hartshorne-Speiser by proving that if $R$ is any regular ring  containing a field of characteristic
zero and $Y\subseteq \Spec R$ is a locally closed subscheme, then ${\rm inj}\dim H_Y^i(R)\leq \dim\Supp H_Y^i(R)$.

As a final main result, we able to  obtain a lower bound for the injective dimension of the local cohomology module $H^{\height I}_I(R)$, in the case
$(R, \frak m)$ is a complete  equidimensional local ring.  More precisely, we show that:

\begin{thm}
Let $(R, \m)$ be a complete local equidimensional ring and $I$ an ideal of $R$.  Then
${\rm inj}\dim H_I^{\height I}(R)\geq \dim R-\height I.$ In particular, if  $R$ is a regular local ring containing a field, then
${\rm inj}\dim H_I^{\height I}(R)=\dim R-\height I$.
\end{thm}
Finally, we will end the paper with an example, which shows that   Theorem 1.4  does not hold in general. \\

  For each $R$-module $L$, we denote by
 ${\rm Ass h}_R(L)$ (resp. ${\rm mAss}_RL$) the set $\{\frak p\in \Ass
_R(L):\, \dim R/\frak p= \dim L\}$ (resp. the set of minimal primes of
$\Ass_RL$).   Also,  the set of all zerodivisors on $L$ is denoted by $Z_R(L)$. For any ideal $\frak{b}$ of
$R$, {\it the radical of} $\frak{b}$, denoted by $\Rad(\frak{b})$,
is defined to be the set $\{x\in R \,: \, x^n \in \frak{b}$ for some
$n \in \mathbb{N}\}$ and we denote $\{{\frak p} \in {\Spec} (R) :
{\frak p} \supseteq {\frak b}\}$ by $V({\frak b})$. Finally, for any
ideal $\frak{b}$ of $R$, {\it the cohomological dimension} of  an
$R$-module $M$, with respect to $\frak{b}$ is defined as
 $${\rm cd}(\frak{b}, M):=\sup \{i\in \mathbb{Z} : H^{i}_{\frak{b}}(M)\neq 0 \}.$$

For any unexplained notation and terminology we refer the reader
to \cite{BS} and \cite{Mat}.

\section{The Results}

The following lemmas  will be quite useful in the proof of the main results. Following $D:=\Hom_R(\bullet, E_R(R/\m))$ (resp.  $\omega_R$) denotes the Matlis duality functor (resp. the canonical module for $R$) (see \cite[3.3]{BH}.
\begin{lem}
 \label{2.1}
Let $(R,\mathfrak{m})$ be a local Noetherian ring and $M$ a finitely
generated $R$-module. Let $\frak p$ be a prime ideal of $R$ such
that $\dim R/{\frak p}=1$ and let $t\geq1$ be an integer. Then
$H^{t}_{\frak m}(M)$ is ${\frak p}$-cofinite if and only if $(H^{t-1}_{\frak p}(M))_{\frak p}=0$.
 \end{lem}
\proof  See \cite[Lemma 2.1]{BAB}. \qed\\

\begin{lem}
Let $(R,\mathfrak{m})$ be a  Cohen-Macaulay local ring  of dimension $d$. Then the $R$-module $H^d_{\m}(R)$ is indecomposable.
\end{lem}
\proof  Without loss of generality,  we may assume that $R$ is a complete  Cohen-Macaulay local ring. Now, we suppose the contrary and we look for a contradiction. Let $H^d_{\m}(R)=A\oplus B$, where $A$ and $B$ are two non-zero Artinian $R$-modules. Then we have $\omega_R\cong D(A)\oplus D(B)$, where $\omega_R$ denotes the canonical module of $R$. So as the $R$-module  $\omega_R$ is indecomposable, it follows that $D(A)=0$ or $D(B)=0$. Hence $A=0$ or $B=0$, which is a contradiction. \qed\\

 The following result will be useful in the proof of the main results in this section.

\begin{thm}
 Let $(R,\mathfrak{m})$ be a  Cohen-Macaulay local ring and let $X$ and $Y$ be non-empty subsets of $\Ass_R(R)$ such that $\Ass_R(R)=X\cup Y$ and $X\cap Y=\emptyset$. Set $$I:=\cap_{\p\in X}\p\,\,\,\,\,{\rm and}\,\,\,\,\,J:=\cap_{\q \in Y}\q.$$ Then  $R/(I+J)$ is an equidimensional local ring of dimension $\dim R-1$.
\end{thm}

\proof  It follows from the hypothesis $X\cap Y=\emptyset$ that $\height(I+J)\geq 1$. Now, we show that  $\height(I+J)= 1$. To do this, suppose the contrary is true. Then
there exists a minimal prime ideal $\frak p$ over $I+J$ such that $\height  \frak p :=n> 1$.
 Since  $\Ass_R(R)=X\cup Y$ it follows that $I\cap J={\rm nil}(R)$, and so $I\cap J$ is a nilpotent ideal of $R$. Therefore
$$H^{n-1}_{I\cap J}(R)=0 = H^n_{I\cap J}(R).$$
Now, in view of the Mayer-Vietoris sequence (see e.g., \cite[Theorem 3.2.3]{BS}) we obtain the isomorphism $$H^n_{I+ J}(R) \cong H^n_I(R) \oplus H^n_J(R).$$
Therefore
$$H^n_{\frak pR_{\frak p}}(R_{\frak p})= H^n_{(I+J)R_{\frak p}}(R_{\frak p})\cong H^n_{IR_{\frak p}}(R_{\frak p}) \oplus H^n_{JR_{\frak p}}(R_{\frak p}).$$

Now, using Lemma 2.2, we deduce that
 $$H^n_{{\frak p}R_{\frak p}}(R_{\frak p})\cong H^n_{IR_{\frak p}}(R_{\frak p})\,\,\,\,\,{\rm or}\,\,\,\,\, H^n_{{\frak p}R_{\frak p}}(R_{\frak p})\cong H^n_{JR_{\frak p}}(R_{\frak p}). $$
Consequently, in view of \cite[Proposition 5.1]{Me},  $H^n_{{\frak p}R_{\frak p}}(R_{\frak p})$ is an $IR_{\frak p}$ or
$JR_{\frak p}$-cofinite  $R_{\frak p}$-module.  Next, as $\height  \frak p > 1$,  it is easy to see that there exists a prime ideal $\frak q \in V(I)$ or
$\frak q \in V(J)$ such that $\frak q \subseteq \frak p$ and $\height  \frak p /\frak q= 1$. Now, using \cite[Proposition 4.1]{Me}, one easily sees that the
$R_{\frak p}$-module $H^n_{{\frak p}R_{\frak p}}(R_{\frak p})$ is ${\frak q}R_{\frak p}$-cofinite. Therefore, it follows from Lemma 2.1 that $H^{n-1}_{{\frak q}R_{\frak q}}(R_{\frak q})=0$. On the other hand, as $R$ is catenary, it follows that
$\height  \frak p /\frak q= \height   \p - \height  \frak q$, and so $$\height  \frak q= \height  \frak p -1= n-1.$$  Hence in view of  Grothendieck's non-vanishing theorem we have $H^{n-1}_{{\frak q}R_{\frak q}}(R_{\frak q})\neq 0$, which is a contradiction. Therefore $\height  \frak p = 1$, and so
$\height(I+J)= 1$.  Now, as $R$ is Cohen-Macaulay, it follows easily that $R/(I+J)$ is an equidimensional ring of dimension $\dim R-1$, as required.    \qed\\

\begin{cor} Let $(R,\mathfrak{m})$ be a  Cohen-Macaulay local ring  and let $x_1,\dots, x_t$ be an $R$-regular sequence. Let $X$ and $Y$
 be non-empty subsets of $\Ass_R(R/(x_1,\dots, x_t))$ such that $\Ass_R(R/(x_1,\dots, x_t))=X\cup Y$ and $X\cap Y=\emptyset$. Set
$$I:=\cap_{\p\in X}\p\,\,\,\,\,{\rm and}\,\,\,\,\,J:=\cap_{\q \in Y}\q.$$Then  $R/(I+J)$ is an equidimensional local ring of dimension $\dim R-t-1$.
\end{cor}
\proof
Since $R/(x_1,\dots, x_t)$ is a Cohen-Macaulay local ring, the assertion follows easily from Theorem 2.3.\qed\\
\begin{lem}
Let $R$ be a  Noetherian ring and $I$ an ideal of $R$ such that ${\rm cd}(I, R)=n>0$. Then the $R$-module  $H^n_I(R)$ is not finitely generated.
\end{lem}
\proof Since by the definition we have $H^n_I(R)\neq 0$, it follows that $\Supp H^n_I(R)\neq \emptyset$. Let $\p\in \Supp H^n_I(R)$. Then it is easy to see that ${\rm cd}(IR_{\p},R_{\p})=n>0$. So replacing of the ring $R$ with the local ring $(R_{\p},\p R_{\p})$, we may assume that $(R,\m)$ is a Noetherian local ring and $I$ is an ideal of $R$ such that ${\rm cd}(I,R)=n>0$. Then using \cite[Exercise 6.1.8]{BS} and Grothendieck's vanishing theorem we have: $$H^n_I(R)/\m H^n_I(R)\cong H^n_I(R)\otimes_R R/\m\cong H^n_I(R/\m)=0.$$Therefore, $H^n_I(R)=\m H^n_I(R)$ and hence using  Nakayama's lemma we can deduce that the $R$-module $H^n_I(R)$ is not finitely generated.\qed\\

We are now in a position to state and prove the first main result of this paper, which  investigates  the finiteness
dimension $f_I(R)$ over a Cohen-Macaulay local ring.

\begin{thm} Let $(R,\mathfrak{m})$ be a   Cohen-Macaulay local ring and  $I$  a  non-nilpotent ideal of $R$. Then
 $$f_I(R)={\rm max}\{1,\height I \}.$$
\end{thm}
\proof To prove there are two cases to consider:\\

{\bf Case 1.}  Suppose that $\height I=0$. Put
\begin{center}
$X:=\Ass_R(R)\cap V(I)$ and $Y:=\Ass_R(R)\backslash V(I)$.
\end{center} Let
$J:= \cap_{\p\in X}\p$ and $K:=\cap_{\q\in Y}\q$.  Since $I$ is not nilpotent it follows that $Y\neq \emptyset$. Also,  as $\height I=0$, it follows that $X \neq \emptyset$. Moreover,  it is easy to see that $\Ass_R(R)=X\cup Y$.  Hence, in view of the proof of Theorem 2.3, we have $\height (J+K)=1$. Therefore, there exists a minimal prime ideal $\frak p$ over $J+K$ such that $\height \frak p=1$. Since $K\subseteq \p$, there exists an ideal $\q \in Y$ such that $\q \subseteq \p$. As $I\subseteq J \subseteq \p$, it follows that $I+\q \subseteq \p$. Moreover, as $I \not\subseteq \q$ it follows that $\height(I+\q)>0$. Therefore, $\height(I+\q)=\height \p=1$. Thus, $\p$ is a minimal prime ideal over $I+\q$ and so $IR_{\p}+\q R_{\p}$ is a $\frak p R_{\frak p}$-primary ideal. Hence, by  Grothendieck's non-vanishing theorem we have $H^1_{IR_{\p}}(R_{\p}/\q R_{\p})\neq 0$. Consequently, it follows from  Grothendieck's vanishing  theorem that ${\rm cd}(IR_{\p}, R_{\p}/\q R_{\p})=1$.  Now, as $\Supp(R_{\p}/\q R_{\p})\subseteq \Spec R_{\p}$, it follows from \cite[Theorem 2.2]{DNT} that
 $${\rm cd}(IR_{\p}, R_{\p})\geq {\rm cd}(IR_{\p}, R_{\p}/\q R_{\p})=1.$$  By using   Grothendieck's vanishing theorem we can deduce that
${\rm cd}(IR_{\p}, R_{\p})=1$ and so by Lemma 2.5 the $R_{\p}$-module $H^1_{IR_{\p}}(R_{\p})\cong (H^1_I(R))_{\p}$ is not finitely generated. In particular, the $R$-module $H^1_I(R)$ is not finitely generated. Now, as the $R$-module $H^0_I(R)$ is finitely generated, it follows that $$f_I(R)=1={\rm max}\{1,0\}={\rm max}\{1,\height I\},$$  as required.

{\bf Case 2.} Now suppose that $\height I=n\geq 1$. Then we have $\grade I=n$ and so in view of \cite[Theorem 6.2.7]{BS},  $f_I(R)\geq n$. Moreover, by the definition there exists a minimal prime ideal $\frak q$ over $I$ such that $\height \q=n$. Hence, in view of  Grothendieck's vanishing and non-vanishing theorems we have $${\rm cd}(IR_{\q},R_{\q})={\rm cd}({\q}R_{\q},R_{\q})=n.$$ Thus,  by Lemma 2.5,  the $R_{\q}$-module $H^n_{IR_{\q}}(R_{\q})\cong (H^n_I(R))_{\q}$ is not finitely generated. In particular, the $R$-module $H^n_I(R)$ is not finitely generated. Hence in view of the definition we have $$f_I(R)=n={\rm max}\{1,n\}={\rm max}\{1,\height I\},$$ and this completes the proof.\qed\\

The next theorem is the second main result of this paper.

\begin{thm}
Let $(R,\mathfrak{m})$ be a  Cohen-Macaulay local ring  and let $I$ be a non-nilpotent ideal of $R$ such that
$\mAss_R(R/I)\subseteq \Ass_R(R)$. Then  $R/(I+\Gamma_I(R))$ is an equidimensional local ring of dimension $\dim R-1$.
\end{thm}
\proof
Since $I$ is not nilpotent, it is clear that $\Gamma_I(R)\subseteq Z_R(R)$ and so it follows  from \cite[Theorem 17.4]{Mat} that $\dim R/\Gamma_I(R)=\dim R$. Moreover, as $$\Ass_R(R/\Gamma_I(R))=\Ass_R(R) \setminus V(I),$$ it follows that $I$ contains an $R/\Gamma_I(R)$-regular element $x$, and so
$$\dim R/(xR+\Gamma_I(R))=\dim R/\Gamma_I(R)-1.$$
Hence $\dim R/(I+\Gamma_I(R))\leq d-1$.

Next, in view of  the Artin-Rees lemma there exists a positive integer $s$ such that $I^s\cap \Gamma_I(R)=0$ and so
$$H^{n-1}_{I^s\cap \Gamma_I(R)}(R)=0 = H^n_{I^s\cap \Gamma_I(R)}(R).$$
 Hence, the Mayer-Vietoris sequence (see e.g., \cite[Theorem 3.2.3]{BS}) yields the isomorphism
$$H^n_{I+ \Gamma_I(R)}(R)=H^n_{I^s+ \Gamma_I(R)}(R) \cong H^n_{I^s}(R) \oplus H^n_{\Gamma_I(R)}(R) \cong H^n_{I}(R) \oplus H^n_{\Gamma_I(R)}(R).$$

Now,  suppose that $\p$ is a minimal prime ideal over $I+\Gamma_I(R)$ such that $\height \p=n> 1$.
Then, as $\p$ is minimal over $I+\Gamma_I(R)$ we get the following isomorphism
$$H^n_{\p R_{\p}}(R_{\p})\cong H^n_{IR_{\p}}(R_{\p}) \oplus H^n_{\Gamma_{IR_{\p}}}(R_{\p}).$$

Now, using Lemma 2.2, we deduce that
 $$H^n_{\p R_{\p}}(R_{\p})\cong H^n_{IR_{\p}}(R_{\p})\,\,\,\,\,{\rm or}\,\,\,\,\,H^n_{\p R_{\p}}(R_{\p})\cong H^n_{\Gamma_{IR_{\p}}(R_{\p})}(R_{\p}).$$

Assume that $H^n_{{\p}R_{\p}}(R_{\p})\cong H^n_{IR_{\p}}(R_{\p})$. Then, in view of \cite[Proposition 5.1]{Me},  $H^n_{{\frak p}R_{\frak p}}(R_{\frak p})$ is an $IR_{\frak p}$-cofinite  $R_{\frak p}$-module.  Next, as $\height  \frak p > 1$ and $\mAss_R(R/I)\subseteq \Ass_R(R)$,  it is easy to see that there exists a prime ideal $\frak q \in V(I)$ such that $\frak q \subseteq \frak p$ and $\height  \frak p /\frak q= 1$. Now, using \cite[Proposition 4.1]{Me}, it follows easily  that the  $R_{\frak p}$-module $H^n_{{\frak p}R_{\frak p}}(R_{\frak p})$ is ${\frak q}R_{\frak p}$-cofinite. Therefore, it follows from Lemma 2.1 that $H^{n-1}_{{\frak q}R_{\frak q}}(R_{\frak q})=0$. On the other hand, as $R$ is catenary, it follows that
$\height  \frak p /\frak q= \height  \frak p - \height  \frak q$, and so $$\height  \frak q= \height  \frak p -1= n-1.$$  Hence in view of  Grothendieck's non-vanishing theorem we have $H^{n-1}_{{\frak q}R_{\frak q}}(R_{\frak q})\neq 0$, which is a contradiction.

 Now, assume that $H^n_{\p R_{\p}}(R_{\p})\cong H^n_{\Gamma_{IR_{\p}}(R_{\p})}(R_{\p})$. Then, again using the fact that
$$\Ass_R(R/\Gamma_I(R))=\Ass_R(R)\backslash V(I)\subseteq \Ass_R(R),$$
 and repeating the above argument we derive a contradiction.   Therefore $\height  \frak p = 1$, and so
$\height(I+\Gamma_I(R))= 1$.
 Now, as $R$ is Cohen-Macaulay, it follows easily that $R/(I+\Gamma_I(R))$ is an equidimensional local ring of dimension $\dim R-1$, as required.    \qed\\

\begin{cor}
 Let $(R,\mathfrak{m})$ be a  Cohen-Macaulay local ring and let $A$ be a non-empty proper subset of  $\Ass_R(R)$.Then $R/(I+\Gamma_I(R))$
is an  equidimensional local ring of dimension $\dim R-1$, where $I=\cap_{\p\in A}\p$.
\end{cor}
\proof The assertion follows easily from Theorem 2.7.\qed\\

\begin{prop}
 Let $(R,\m)$ be a  Cohen-Macaulay local ring and let  $\Ass_R(R)=\{\p_1,\dots,\p_n\}$, $n\geq 2$.  Let  $A_j=\Ass_R(R)\backslash \{\p_j\}$ and $I_j=\cap_{\p\in A_j}\p$, for all $1\leq j \leq n$. Then $0=\cap_{j=1}^n\Gamma_{I_j}(R)$ is the unique reduced primary decomposition of the zero ideal $0$ in $R$,  $\Gamma_{I_j}(R)$ is a $\frak p_j$-primary ideal of $R$  and $R/(I_j+\Gamma_{I_j}(R))$ is an equidimensional local ring of dimension $\dim R-1$.
\end{prop}
\proof
 As $$\Ass_R(R/\Gamma_{I_j}(R))=\Ass_R(R)\backslash V(I_j)=\{\p_j\},$$ it follows that $\Gamma_{I_j}(R)$ is a $\p_j$-primary ideal of $R$. Now,
 we show that $\cap_{j=1}^n\Gamma_{I_j}(R)=0$. To this end, we assume  that $\cap_{j=1}^n\Gamma_{I_j}(R)\neq 0$ and derive a contradiction.
Let $a\in \cap_{j=1}^n\Gamma_{I_j}(R)$ be such that $a\neq 0$.  Then  $(0:_Ra)\subseteq Z_R(R)$, and so there exists $\p_j\in \Ass_R(R)$ such that
 $(0:_Ra) \subseteq \p_j$. Next,  as $a\in\Gamma_{I_j}(R)$ it follows that  there exists a positive integer $k$ such that
 $I_j^k\subseteq (0:_Ra)\subseteq \p_j$, and  so $I_j\subseteq \p_j$. Therefore, there exists $\p_i\in A_j$ such that $\p_i\subsetneqq \p_j$,
 which is a contradiction, (note that  $\Ass_R(R)= \mAss_R(R)$).  Now, using \cite[Theorem 6.8]{Mat} we see that  $\p_j$-primary component
 $\Gamma_{I_j}(R)$ of the zero ideal $0$ of $R$ is uniquely determined. That is,
 $0=\cap_{j=1}^n\Gamma_{I_j}(R)$ is the unique reduced primary decomposition of the zero ideal $0$ in $R$. Moreover, it follows from  Corollary 2.8 that the ring $R/(I_j+\Gamma_{I_j}(R))$ is  equidimensional local of dimension $\dim R-1$. \qed\\

 The following lemma is needed in the proof of Theorem 2.11.

 \begin{lem}
Let $(R, \m)$ be a  local  ring  and $M$ an arbitrary  $R$-module. Let $x$ be an element of $\m$  such that
$x\not\in \bigcup_{\p \in \Ass_R (M)\backslash V(\m)} \p$. Then $\Gamma_{Rx}(M)=\Gamma_{\m}(M)$.
\end{lem}
 \proof As $x\in \m$, it is enough to show that $\Gamma_{Rx}(M)\subseteq \Gamma_{\m}(M)$. To do this, let $w\in \Gamma_{Rx}(M)$.
 Then $x\in \Rad(0:_Rw)$.  Since $\mAss_RR/(0:_Rw)\subseteq \Ass_R(M)$, it follows from the assumption
  $x\not\in \bigcup_{\p \in \Ass_R (M)\backslash V(\m)} \p$
 that $\Rad(0:_Rw)=\m$, and so there exists $n\in \mathbb{N}$ such that $\m^nw=0$. Thus $w\in \Gamma_{\m}(M)$, as required. \qed\\

The following theorem is in preparation for the third main result of this paper, which gives us a lower bound of injective dimension of $ H_I^{\height I}(R)$.
  Here $D_I(R)$
denotes the ideal transform of $R$ with respect to $I$ (see \cite[2.2.1]{BS}).

\begin{thm}
Let $(R, \m)$ be a complete local equidimensional ring of dimension $d$ and $I$ an ideal of $R$ such that $\height I=t$. Then
$H_{\m}^{d-t}(H_I^t(R))\neq0$.
In particular, $${\rm inj}\dim H_I^t(R)\geq d-t.$$
\end{thm}
\proof As $R$ is catenary, it follows from \cite[ Lemma  2, P. 250]{Mat} that
 $$\height J+ \dim R/J=\dim R,$$
for every ideal $J$ of $R$.  In particular,  we have $\dim R/I=d-t$. We now use induction on $d-t$. When $d=t$, the ring $R/I$ is Artinian and so $\Rad(I)=\m$.
Hence $H_I^t(R)=H_{\m}^t(R)$ and so as $H^0_{\m}(H_I^t(R))= H_{\m}^t(R)$, the assertion follows from  Grothendieck's  non-vanishing theorem (see  \cite[Theorem 6.1.4]{BS}) in this case.

Assume, inductively, that $d-t>0$ and that the result has been proved for the ideals $J$ with $\dim R/J=0,1, \dots, d-t-1$. Since the sets $\Ass_R(H^t_I(R))$ and
 $\Ass_R (H^{t+1}_I(R))$ are countable, it follows from \cite[Lemma 3.2]{MV} that
 $$\m\nsubseteq (\bigcup_{\p \in \Ass_R(H^t_I(R))\backslash V(\m)} \p)\, \bigcup\, (\bigcup_{\p \in \Ass_R(H^{t+1}_I(R))\backslash V(\m)} \p)\,\bigcup\, (\bigcup_{\p \in {\rm Assh}_R(R/I)} \p).$$
 Whence, there exists $x\in \frak m$ such that
 $$x\not\in  (\bigcup_{\p \in \Ass_R(H^t_I(R))\backslash V(\m)} \p)\, \bigcup\, (\bigcup_{\p \in \Ass_R(H^{t+1}_I(R))\backslash V(\m)} \p)\,\bigcup\, (\bigcup_{\p \in {\rm Assh}_R(R/I)} \p).$$
 Then it follows easily from $x\not\in \bigcup_{\p \in {\rm Assh}_R(R/I)} \p$ that $$\dim R/(I+Rx)=d-t-1,$$
  and in view of Lemma 2.10 we have
 \begin{center}
 $\Gamma_{Rx}(H_I^t(R))=\Gamma_{\m}(H_I^t(R))$ \,\,\,\,\, and \,\,\,\,\, $\Gamma_{Rx}(H_I^{t+1}(R))=\Gamma_{\m}(H_I^{t+1}(R))$.
 \end{center}

Moreover, there is an exact sequence

$$0\longrightarrow H^{1}_{Rx}(H^{t}_{I}(R))\longrightarrow H^{t+1}_{I+Rx}(R) \longrightarrow H^{0}_{Rx}(H^{t+1}_{I}(R))\longrightarrow 0, \hspace{13mm}(\dag)$$
(see \cite[Corollary 3.5]{Sc}).

Now, if $\dim R/I=1$ then  in view of  \cite[Theorem 2.6]{BN1} the $R$-module $H^{t}_{I}(R)=H^{d-1}_{I}(R)$ is $I$-cofinite.  Next, it is easy
to see that  $\dim\Supp H_I^{d-1}(R)=1$, note that $\dim R/I=1$. Hence, it follows from \cite[Theorem 2.9]{Ma} that  $H^1_{\m}(H_I^{d-1}(R))\neq 0$,
and so the result has been proved in this case. Therefore, we assume that $\dim R/I\geq2$. Then $$\dim R/(I+Rx)=d-t-1\geq1,$$ and so in view of  Grothendieck's vanishing theorem $$H_{\m}^{d-t-1}(\Gamma_{\m}(H_I^{t+1}(R)))=0.$$
Hence by using the exact sequence $(\dag)$ we obtain the following exact sequence
$$H_{\m}^{d-t-1}(H^{1}_{Rx}(H^{t}_{I}(R)))\longrightarrow H_{\m}^{d-t-1}(H^{t+1}_{I+Rx}(R) )\longrightarrow 0.$$
Thus by the inductive hypothesis $H_{\m}^{d-t-1}(H^{1}_{Rx}(H^{t}_{I}(R)))\neq0$.

On the other hand, since $d-t>0$, it yields that
$$H_{\m}^{d-t}(H^{t}_{I}(R))\cong H_{\m}^{d-t}(H^{t}_{I}(R)/\Gamma_{\m}(H^{t}_{I}(R))).$$
Now, let $T:=H^{t}_{I}(R)/\Gamma_{\m}(H^{t}_{I}(R))$. It is thus sufficient for us to show that  $H^{d-t}_{\m}(T)\neq0$.
To do this, in view of \cite[Remark 2.2.17]{BS}, there is the exact sequence
$$0 \longrightarrow T \longrightarrow D_{Rx}(T) \longrightarrow H^{1}_{Rx}(T)\longrightarrow 0. \hspace{13mm}(\dag\dag)$$
Also, in view of \cite[Theorem 2.2.16]{BS}, we have $ D_{Rx}(T)\cong T_x$, and so $$ D_{Rx}(T)\stackrel{x}\longrightarrow  D_{Rx}(T),$$
is an $R$-isomorphism. Therefore, for all $i\geq0$, $$ H_{\m}^i(D_{Rx}(T))\stackrel{x}\longrightarrow H_{\m}^i(D_{Rx}(T)),$$
 is an $R$-isomorphism, and hence $ H_{\m}^i(D_{Rx}(T))=0$, for all $i\geq0$. Consequently, it follows from the exact sequence $(\dag\dag)$ that
$$H^{d-t}_{\m}(T)\cong H^{d-t-1}_{\m}(H^1_{Rx}(T)).$$
As $H^{d-t-1}_{\m}(H^1_{Rx}(T))\neq0$,  this completes the inductive step. \qed\\

\begin{cor} Let $(R,\mathfrak{m})$ be a  Cohen-Macaulay local ring of dimension $d$ and $I$ an ideal of $R$ such that $\height I=t$. Then
$H_{\m}^{d-t}(H_I^t(R))\neq0$. In particular, ${\rm inj}\dim H_I^t(R)\geq d-t$.
\end{cor}
\proof Let $\hat{R}$ denote the completion of $R$ with respect to the $\frak m$-adic topology. Then, as $(\hat{R},\mathfrak{m}\hat{R})$ is a  complete local equidimensional ring of dimension $d$, the assertion follows from Theorem 2.11, the faithfully flatness of the homomorphism $R\longrightarrow \hat{R}$ and the fact that  $$\height I=\grade I= {\rm grade}\, I \hat{R}=\height I\hat{R}.$$ \qed

\begin{lem}
Let $(R, \m)$ be a regular local ring containing a field and $I$ an ideal of $R$. Then, for any integer $n$ with $H^n_I(R)\neq0$,
$${\rm inj}\dim H_I^n(R)\leq \dim\Supp H_I^n(R).$$
\end{lem}
\proof The result follows from \cite{HS} and \cite{Ly1}.  \qed\\

\begin{cor}
Let $(R, \m)$ be a regular local ring containing a field and $I$ an ideal of $R$ such that $\height I=t$. Then
$${\rm inj}\dim H_I^t(R)=\dim R-t.$$
\end{cor}
\proof In view of Corollary 2.12 and Lemma 2.13, it is enough to show that $$\dim\Supp H_I^t(R)=\dim R-t.$$
To this end, as $\Supp H_I^t(R)\subseteq V(I)$ and $\dim R/I=\dim R-t$, we have $$\dim\Supp H_I^t(R)\leq\dim R-t.$$
On the other hand, since $\height I=t$ there exists a minimal prime $\p$ over $I$ such that $\height \p=t$. Now, in view of
\cite[Theorems 4.3.2 and 6.1.4]{BS} we deduce that
$$(H_I^t(R))_{\p}\cong H_{IR_{\p}}^t(R_{\p})\cong H_{\p R_{\p}}^t(R_{\p})\neq0.$$
Thus $\p \in \Supp H_I^t(R)$,  and so as $\dim R/\p=\dim R-t$, it follows that   $$\dim\Supp H_I^t(R)\geq\dim R-t.$$
This completes the proof.  \qed\\

We end the paper with the following example, which shows that  Corollary 2.14 does not hold in general.

\begin{exam}
Let $(R, \m)$ be a regular local ring containing a field with $\dim R=d\geq 3$,  $\p$  a prime ideal of $R$ such that $\dim R/\p=1$ and
$x\in \m\backslash \p$. Then ${\rm inj}\dim H^{\dim R-1}_{Rx \cap \frak p }(R)=0,$ and  $\dim\Supp H^{\dim R-1}_{Rx \cap \p }(R)=1.$
\end{exam}
\proof    Since $\Rad(\p+Rx)=\m$, it follows from the Mayer-Vietoris sequence (see e.g., \cite[Theorem 3.2.3]{BS}) that
$$0\longrightarrow H^{d-1}_{\p}(R)\longrightarrow H^{d-1}_{x\p}(R) \longrightarrow H^d_{\m}(R)\,\,\,\,\,\,\,\,\,\,\,\,\,\,\,\, (\dag\dag\dag)$$
is an exact sequence. Since, in view of the proof of Corollary 2.14, $\dim\Supp H^{d-1}_{\p}(R)=1$ and $H^d_{\m}(R)$ is Artinian, it follows that $\dim\Supp H^{d-1}_{x\p}(R)=1$.  

On the other hand,  the exact sequence
$$0\longrightarrow R  \stackrel{x} \longrightarrow R \longrightarrow R/xR \longrightarrow 0$$
 induces the exact sequence
$$ H^{d-2}_{x\p}(R/xR) \longrightarrow H^{d-1}_{x\p}(R) \stackrel{x} \longrightarrow H^{d-1}_{x\p}(R) \longrightarrow H^{d-1}_{x\p}(R/xR).$$

Since $\Gamma_{ x\p}(R/xR)=R/xR$ and $d\geq3$, it follows that
$$H^{d-2}_{ x\p}(R/xR)=0=H^{d-1}_{ x\p}(R/xR).$$
Therefore,  the $R$-homomorphism $$H^{d-1}_{x\p}(R)\stackrel{x} \longrightarrow H^{d-1}_{x\p}(R)$$
is an isomorphism, and so $(H^{d-1}_{x\p}(R))_x\cong H^{d-1}_{x\p}(R)$.

 On the other hand, from the exact sequence $ (\dag\dag\dag)$, we have $$(H^{d-1}_{x\p}(R))_x\cong (H^{d-1}_{\p}(R))_x.$$

 Moreover, the exact sequence $$0 \longrightarrow H^{d-1}_{\p }(R) \longrightarrow E_R(R/\p) \longrightarrow E_R(R/\m),  $$
implies that $$(H^{d-1}_{\p}(R))_x \cong (E_R(R/\p))_x \cong E_R(R/\p).$$ Therefore $H^{d-1}_{x\p}(R) \cong E_R(R/\p)$,  and so
${\rm inj}\dim H^{d-1}_{x\p}(R)=0$, as required.  \qed\\

\begin{center}
{\bf Acknowledgments}
\end{center}
The authors are deeply grateful to the referee for his/her careful reading and  helpful suggestions on
the paper.  We also would like to thank Professor Hossein Zakeri for his reading of the first draft and valuable discussions.
 Finally, the authors would like to thank from the Institute for Research in Fundamental Sciences (IPM) for the financial support.

\end{document}